\documentclass[12pt]{article}

\font\tenmsbm=msbm10 scaled 1200
\font\sevenmsbm=msbm9
\newfam\msbmfam
\textfont\msbmfam=\tenmsbm
\scriptfont\msbmfam=\sevenmsbm
\def\msbm{\fam\msbmfam\tenmsbm}

\def\beq{\begin{equation}}
\def\eeq{\end{equation}}
\def\a{\alpha}

\def\O{\cal O}

\def\bc{\begin{center}}
\def\ec{\end{center}}

\def\vuoto{\ \hfill\hbox{\vbox{\hrule\hbox{\vrule
height5pt\kern5pt\vrule height5pt}\hrule}}\par\medskip\rm}  

\def\P{\mbox{\msbm P}}
\def\r{\longrightarrow}

\begin{document}

\title{On two simple criteria for recognizing complete intersections in 
codimension 2}

\author{ Alessandro Arsie \thanks{e-mail: arsie@sissa.it} \\
S.I.S.S.A. - I.S.A.S. \\
Via Beirut, 4 - 34013 Trieste,  \\
and \\
Dipartimento di Matematica, Universit\'a di Ferrara \\
Via Machiavelli, 35, 44100 Ferrara, Italy}

\date{}
\maketitle

\begin{abstract}
Developing a previous idea of Faltings,
we characterize the complete intersections of codimension 2 in 
${\hbox{\msbm{ 
P}}}^{n}$, $n\geq 3$, over an algebraically closed field of any 
characteristic, among l.c.i. $X$, as those that are subcanonical 
and scheme-theoretically defined by $p\leq n-1$ equations. Moreover, we 
give some other results assuming that the normal bundle of $X$ extends 
to a numerically split bundle on ${\hbox{\msbm{P}}}^{n}$, $p\leq n$ and the characteristic of the base 
field is zero. Finally, using our characterization,
we give a (partial) answer to a question posed recently by Franco, 
Kleiman and Lascu (\cite{FKL}) on self-linking and complete intersections in 
positive characteristic.

\end{abstract}

{\small MSC (1991): Primary: 14M10 , Secondary: 14M06}

{\small Keywords: local complete intersections, complete intersections, 
Serre correspondence, self-linking.}
\bc

\ec

\section{Introduction}

 The problem of detecting (global) complete intersections is a key question in 
 projective algebraic geometry and commutative algebra. Up to now, 
 this problem is far from being solved and a complete answer is known only 
 in trivial cases, such that of hypersurfaces in projective spaces or 
 in Grassmanians. Moreover, the relevant conjecture of Hartshorne 
 according to which projective subvarieties of small codimension, 
 compared to their dimension, have to be complete intersections, is 
 still not proved. However, in the last 25 years, there have been 
 some partial results in this direction, particularly in the case of 
 codimension 2. Essentially, the results obtained in the case $X$ is a 
 smooth codimension 2 subvariety of ${\P}^{n}$, $n\geq 6$ can be grouped 
 into two kinds of criteria.
 
  The first one says that if $X$ is 
 contained in a hypersurface $V$, such that $deg(V)\leq n-2$, then $X$ 
 is a complete intersection (see \cite{R} or the recent improvement 
 in  \cite{EF}, where it is shown that the bound on the degree of $V$ 
 can be increased to $n-1$, in the case of codimension 2 
 subvarieties of ${\P}^{6}$); using this kind of 
 criterion one can give also a bound on the degree of $X$, so that to 
 assure that $X$ is a complete intersection.    
 
 The second kind of criterion is based on giving a bound on the 
 number $p$ of generators, not for the homogeneous ideal $I(X)$, but for an ideal 
 $I_{sch}(X)$ which coincide with $I(X)$ only in high degree, that is 
 $[I_{sch}(X)]_{d}=[I(X)]_{d}$, for $d \gg 0$. We call $I_{sch}(X)$ 
 the ÒschematicÓ ideal of $X$, in that its generators define $X$ 
 scheme-theoretically. Following this approach, Faltings proved in 
 \cite{F} that if $p\leq n-2$ and $n\geq 8$ and $X$ is a (possibly 
 singular) subcanonical local complete intersection,  then it is a complete 
 intersection (in any characteristic). Some years later, this result was 
 improved in \cite{N}, proving in characteristic zero 
 that if $p\leq n-1$, $n\geq 8$, then $X$ is a complete intersection 
 (but assuming $X$ smooth). 

The aim of the present work is twofold: on one hand we would like to 
give a different (in that we use Serre's correspondence) and simpler proof of the result announced in 
\cite{N}, hoping to give a Òcrystal clearÓ version of some obscure 
(to our opinion) arguments. Moreover, assuming only that $X$ is a 
(possibly singular) subcanonical l.c.i.,  we prove in any 
characteristic that if $n\geq 3$ and $p\leq n-1$, then $X$ is a complete
intersection and we give some more results working in characteristic 
zero, assuming that the normal bundle of $X$ extends to a numerically 
split bundle $E$ on ${\hbox{\msbm{P}}}^{n}$ (i.e. the Chern classes 
of $E$ are those of a split bundle), $n\geq 3$,
 $p\leq n$. On the other hand, as an application of 
our result we answer to a question posed recently by Franco, Kleiman 
and Lascu in \cite{FKL}, (neglecting the case of space curves).
Unfortunately, our result shows that the characterization given by 
Faltings {\em is not peculiar} of two codimension embeddings in high 
dimensional projective spaces.   
  
  Our proof is based in constructing and exploiting an exact 
  sequence of locally free sheaves, (sequence (\ref{fund})), which 
  relates the rank 2 vector bundle $E$ appearing in Serre's 
  correspondence with the generators of the Òscheme-theoreticÓ ideal of $X$.

\section{Main result: a scheme-theoretic criterion}
From now on, $X$ will denote a codimension 2 {\em 
subcanonical} l.c.i. (possibly singular) closed subscheme of a 
projective space ${\msbm P}^{n}$ over an algebraically 
closed field $k$ of any characteristic $p\geq 0$, where, as usual ${\msbm P}^{n}=
Proj (k[x_{0},\ldots,x_{n}])$.
Then we prove the following result:

{\bf Theorem}: {\em A) If $X \subset {\msbm P}^{n}$, ($n\geq 3$) is 
a scheme-theoretic intersection of $p\leq n-1$ hypersurfaces, then $X$ is a complete 
intersection. 

B) Assume moreover that char($k$)=$0$, $n\geq 3$, and $X$ 
scheme-theoretically defined 
by $p\leq n$ hypersurfaces $V_{1},\ldots,V_{n}$ of degrees 
$d_{1},\ldots, d_{p}$, respectively. If the normal bundle of $X$ can 
be extended to a rank 2 vector bundle $E$ on ${\hbox{\msbm{P}}}^{n}$ 
which is numerically split (i.e. $c_{1}(E)=a+b$ and $c_{2}(E)=ab$, $a, 
b\in {\hbox{\msbm {Z}}}$) and $a$} or {\em $b$ is in $(d_{1},\ldots,d_{p})$,
 then $X$ is a complete intersection.}

Proof of part A: Since $X$ is assumed to be subcanonical (i.e. its 
dualizing sheaf $\omega_{X}$, which is locally free, 
is of the form ${\cal O}_{X}(e)$) , by Serre's 
correspondence there exists 
an algebraic vector bundle $E$ of rank 2  over ${\msbm P}^{n}$ and a 
section $s\in H^{0}(\P, E)$ such that $X$ is identified with the 
scheme of zeroes of $s$, $Z(s)$. The Koszul complex for this section 
gives a projective resolution of the ideal sheaf of $Z(s)$, hence of 
the ideal sheaf of $X$:
\beq
\label{koszul1}
0  \r  \bigwedge^{2} E^{*}  \r 
E^{*} 
\r  {\cal I}_{X}  \r  0.
\eeq
Since ${\cal I}_{X}$ is not itself projective, by (\ref{koszul1}) it 
turns out that the projective dimension of ${\cal I}_{X}$ is 1. On 
the other hand, if $X$ is schematically cut out by $p\leq n-1$ 
hypersurfaces of degrees $d_{1},\ldots,d_{p}$, we have an exact sequence:
\beq
\label{koszul2}
0  \r  Ker(f)  \r  \bigoplus {\O} (-d_{i}) \stackrel{f}{\r} 
{\cal I}_{X}  \r   0.
\eeq
Since pd(${\cal I}_{X}$)=1, then the first syzygy $Ker(f)$ is also 
projective (see for example \cite{W}), hence it corresponds to a locally free sheaf. 
Certainly, we have a morphism $h\in Hom(\oplus {\O} (-d_{i}) \oplus 
{\O} (-c_{1}), {\cal I}_{X})$ which is given first by projecting to 
$\oplus {\O} (-d_{i})$ and then composing with $f$, ($c_{1}$ is the 
first Chern class of $E$). Moreover, since $Ext^{1}(\oplus {\O} (-d_{i}) \oplus 
{\O} (-c_{1}), {\O} (-c_{1}))=0$, then $h$ comes from an element $g$ in 
$Hom(\oplus {\O} (-d_{i}) \oplus {\O} (-c_{1}), E^{*})$; indeed, due to 
the fact that $\bigwedge^{2}E^{*}\equiv {\O} (-c_{1})$, we have the 
following commutative diagram:
\[
\begin{array}{ccccccccc}
0 & \r & {\O} (-c_{1}) & \r & \bigoplus {\O} (-d_{i})\oplus {\O} (-c_{1}) & 
\r & \bigoplus {\O} (-d_{i}) & \r & 0 \\

 &  & \downarrow \equiv &  &  \downarrow g  &  & \downarrow f &  &  \\
0 & \r & \bigwedge^{2} E^{*} & \r & E^{*} & \r & {\cal I}_{X} & \r 
& 0
\end{array}
\]
Applying  the snake lemma to the previous commutative diagram, we see 
that $g$ is surjectve and that $Ker(g)\equiv Ker(f)$, so that 
dualizing the sequence $0\r Ker(g) \r \bigoplus {\O} (-d_{i})\oplus {\O} 
(-c_{1}) \stackrel{g}{\r} E^{*}\r 0$, we get a short exact sequence 
of {\em locally free} sheaves:
\beq
\label{fund}
0  \r E  \stackrel{\a_{1}}{\r}  \bigoplus {\O} (d_{i})\oplus {\O} (c_{1}) 
\stackrel{\a_{2}}{\r}  C  \r 
  0,
\eeq
where $C$ is just the cokernel sheaf. Since $C$ is locally free, it 
can be identified with a vector bundle of rank equal to $p-1$ ($p\leq 
n-1$). Let $in_{j}$ denote the canonical injection of ${\O}(d_{j})$ 
into  $\bigoplus {\O} (d_{i})\oplus {\O} (c_{1})$, and $pr_{j}$ the 
corresponding projection from $\bigoplus {\O} (d_{i})\oplus {\O} 
(c_{1})$ to ${\O}(d_{j})$. Considering the maps $f_{j}:=pr_{j}\circ\a_{1}$ 
and $g_{j}:=\a_{2}\circ in_{j}$ we have a diagram like the following:
\[
\begin{array}{ccccccccc}
0 & \r & E & \stackrel{\a_{1}}{\r} & \bigoplus {\O} (d_{i})\oplus {\O} 
(c_{1}) & \stackrel{\a_{2}}{\r} & C & \r & 0 \\
 &    &   &   \stackrel{f_{j}}{\searrow} & pr_{j} 
 \downarrow \quad \uparrow in_{j}  & 
 \stackrel{g_{j}}{\nearrow}&  & &  \\
 &    &   &       &  {\O}(d_{j}) &  &   &   &  
 \end{array}
 \]
Now consider the morphisms $E\stackrel{f_{1}}{\r} {\O}(d_{1})$ and
${\O}(d_{1}) \stackrel{g_{1}}{\r} C$ and denote by $Z(f_{1})$ and 
$Z(g_{1})$, their respective degeneracy loci. Since in general we have 
that
$codim(Z(f_{j}))\leq 2$ and $codim(Z(g_{j}))\leq p-1$, it turns out 
that if $p\leq n-1$, then $Z(f_{1})\cap Z(g_{1})\neq \emptyset$. On 
the other hand, by exactness of (\ref{fund}), it is clear that $
Z(f_{1})\cap Z(g_{1})= \emptyset$.
Indeed, if it exists  $x\in Z(f_{1})\cap Z(g_{1})$, then 
$Ker(f_{1})_{x}=E_{x}$, but since the morphism $\alpha_{1}$ can not 
degenerate at any point (due to the fact that the cokernel $C$ is 
locally free), we have that $Im(f_{1})_{x}\subset \oplus_{i\geq 
2}{\O}_{x}(d_{i})\oplus {\O}_{x}(C_{1})$. On the other hand, 
$Ker(g_{1})_{x}={\O}_{x}(d_{1})$, but since $in_{1}$ is always 
injective, we have that $Ker(\alpha_{2})_{x}=in_{1}({\O}(d_{1}))_{x}$. 
Hence, if exists $x\in Z(f_{1})\cap Z(g_{1})$, then 
$Ker(\alpha_{2})_{x}\cap Im(\alpha_{1})_{x}=\emptyset$ and so the 
sequence (\ref{fund}) can not be exact in the middle, at $x$. Absurd.

Hence $Z(f_{1})=\emptyset$ or $Z(g_{1})=\emptyset$. 

If $Z(f_{1})=\emptyset$, then $f_{1}$ is never 
degenerate, so dualizing $E \stackrel{f_{1}}{\r} {\O}(d_{1}) \r 0$ we 
get $0 \r {\O}(d_{1}) \stackrel{(f_{1})^{T}}{\r} E^{*}$, where the   
map $(f_{1})^{T}$ is never degenerate; hence $E^{*}$ splits, so $E$ 
splits and $X$ is a complete intersection.

If instead $Z(g_{1})=\emptyset$, we build up the following 
commutative diagram:
\[
\begin{array}{ccccccccc}
   &   & 0 & \r & {\O}(d_{1}) & \stackrel{\equiv}{\r} & {\O}(d_{1}) & 
 \r & \ldots \\
   &    & \downarrow &  & \downarrow in_{1} &   &  \downarrow 
   g_{1} &  &   \\
 0 & \r & E & \r & \bigoplus {\O}(d_{i}) \oplus {\O}(c_{1}) & \r & 
 C & \r & 0 \\
   &    &  \downarrow &  & \downarrow &         &  \downarrow
   &    &    \\
 0 & \r & Ker(\psi) & \r & \bigoplus_{i\geq 2} {\O}(d_{i})\oplus 
 {\O}(c_{1}) & \stackrel{\psi}{\r} & C^{'} & \r & 0 \\
  &   & \downarrow &   & \downarrow &    & \downarrow &  &   \\
 \ldots & \r & 0 & \r  & 0 & \r  & 0 & &  \\
\end{array}
\]                       
By applying the snake lemma to the two central rows, we see that 
$Ker(\psi)\equiv E$, $C^{'}$ is locally free since 
$Z(g_{1})=\emptyset$, and we obtain a short exact sequence of {\em 
locally free} sheaves:
\beq
\label{fund2}
0  \r  E \r  \bigoplus_{i\geq 2} {\O}(d_{i})\oplus 
{\O}(c_{1})  \r  C^{'}  \r  0.
\eeq
Repeating the previous reasoning, we can consider the morphisms
$E \stackrel{f_{2}}{\r} {\O}(d_{2})$ and ${\O}(d_{2}) 
\stackrel{g_{2}}{\r} C^{'}$ and as before $Z(f_{2})=\emptyset$ or
$Z(g_{2})=\emptyset$. If $Z(f_{2})=\emptyset$, then $E$ splits and 
$X$ is a complete intersection. On the other hand, if 
$Z(g_{2})=\emptyset$, arguing as before, we obtain a short exact 
sequence of locally free sheaves:
\[ 
0  \r  E \r \bigoplus_{i\geq 3} {\O}(d_{i})\oplus 
{\O}(c_{1})  \r  C^{''}  \r   0. \]
In this way we obtain a sequence $Z(f_{1})$, \ldots, $Z(f_{p-1})$. If 
one of these is empty, we are done; otherwise, if all are  not empty,
 then, necessarily, $Z(g_{p-1})=\emptyset$ and as before 
we obtain:
\[ 0  \r  E  \r  {\O}(d_{p})\oplus 
{\O}(c_{1})  \r  0, \]  
and we are done.
\vuoto

Proof of part B: To deal with the more restrictive case of 
characteristic zero, we can assume $k={\msbm C}$ in view of 
Lefschetz's principle. 
In our assumption
$X$ is cut out schematically by $n$ hypersurfaces 
$V_{1},\ldots,V_{n}$ of degrees $d_{1},\ldots,d_{n}$ and we have that 
 $c_{1}(E):=a+b=d_{k}+b$ and 
$c_{2}(E):=ab=d_{k}b$ for some $k\in (1,\ldots,n)$. (It is not 
restrictive to assume that $a\in (d_{1},\ldots,d_{n})$). 
Reordering the hypersurfaces, we can assume that 
$c_{1}(E)=d_{1}+b$ and $c_{2}(E)=d_{1}b$. From the exact 
sequence (\ref{fund}), it is clear that the rank of $C$ is $n-1$, so 
that the morphism ${\O}(d_{1}) \stackrel{g_{1}}{\r} C$ degenerates at 
most in codimension $n-1$.  On the other hand, the morphism $E 
\stackrel{f_{1}}{\r} {\O}(d_{1})$ can not degenerate in codimension 
$2$, otherwise, if $Z(f_{1})$ is its degeneracy locus, the Poincar\'e 
dual $[Z(f_{1})]\in H^{4}(\P, {\msbm 
Z})$ would represent $c_{2}(E^{*}\otimes 
{\O}(d_{1}))=c_{1}(E^{*})d_{1}+c_{2}(E^{*})+d_{1}^{2}=0$, so that 
either the morphism $f_{1}$ does not degenerate at all, and in this case we 
are done as before, or it degenerates in codimension $1$.

So if the morphism $f_{1}$ degenerates in codimension 
one, we have that $Z(f_{1})\cap Z(g_{1})\neq \emptyset$, provided 
that $Z(g_{1})\neq \emptyset$. On the other hand, by exactness of 
(\ref{fund}), we must have $Z(f_{1})\cap Z(g_{1})= \emptyset$, so that 
we conclude that $Z(g_{1})=\emptyset$. Thus, arguing as in part A, we 
obtain the following short exact sequence:
\[ 0  \r  E \r  \bigoplus_{i\geq 
2}^{n}{\O}(d_{i})\oplus{\O}(c_{1})  \r  C^{'}  \r  0, \]
and, from this, we conclude as in part A.   
\vuoto

The result of part B of the previous theorem can be interpreted as a 
relation between degree $deg(X)$ and subcanonicity $e$, 
recalling the well-known fact that  if $E$ is the vector bundle associated to $X$ 
via Serre's correspondence, then $deg(X)=c_{2}(E)$, while 
$e+n+1=c_{1}(E)$. 

{\bf Corollary A}: {\em Let $X\subset {\P}^{n}$ ($X$ as in the 
hypotheses of Theorem B), $n\geq 3$ be 
scheme-theoretically defined by $n$ hypersurfaces of degrees 
$d_{1},\ldots,d_{n}$, 
and let $l$ be an integer in the set $(d_{1},\ldots,d_{n})$. 
If the following relation is satisfied:}
\beq
\label{last}
deg(X)+l^{2}-(e+n+1)l=0,
\eeq
{\em then $X$ is a complete intersection}.

Proof: 
Arguing as in part B, it is clear that to show that $E$ splits is 
sufficient to show that $c_{2}(E^{*}\otimes {\O}(l))=0$ for some $l$ 
as above. But the vanishing of the 
second Chern class of $E^{*}\otimes {\O}(l)$ is exactly the 
relation (\ref{last}), as an easy computation can show. 
\vuoto

{\bf Remark 1}: The approach of giving bounds on the degree of a 
subvariety to detect a complete intersection is particularly 
"effective", but  it is  obviously hopeless if one pretend to solve 
Hartshorne's conjecture. 
On the other hand, since any closed subscheme (irreducible or not) of ${\P}^{n}$ 
which is a local complete intersection is always 
scheme-theoretically defined by $n+1$ hypersurfaces, as proved in 
\cite{Fu}, the approach of recognizing a complete intersection via the 
number of generators of its scheme-theoretic ideal, could be in 
principle useful 
to solve the conjecture. Unfortunately, the cases $p=n$ and in 
particular $p=n+1$ (the generic case) appear completely intractable, 
at least up to now,
since it is very difficult to relate the algebro-geometric properties 
of a small codimension embedding, with those of its scheme-theoretic 
ideal.

{\bf Remark 2}: In the light of the previous remark and of Theorem B, 
it would be nice to know when a subvariety can be scheme-theoretically 
defined by $n$ equations. To get a sufficient condition, we can use 
the theory of excess and residual intersections as developed in 
\cite{F}. For example, let us consider $4$ hypersurfaces 
$\{V_{1},\ldots,V_{4}\}$ in ${\P}^{4}$, 
such that $\cap V_{i}=X\cup \{p_{1},\ldots,p_{k}\}$, where $X$ is a 
smooth subcanonical surface and $\{p_{1},\ldots,p_{k}\}$ are 
(possibly non reduced) points,
i.e. the four hypersurfaces define scheme-theoretically the union of 
$X$ and a bunch of points outside $X$. The theory of residual 
intersections enables us to predict the (weighted) number of residual 
points as a function of the degrees of the hypersurfaces, of the degree 
of $X$ and of the degrees of the Chern classes of $T_{X}$, the 
tangent sheaf of $X$. Imposing that the number of residual points is 
zero, gives us a sufficient condition for a surface (a subvariety in 
general) to be scheme-theoretically defined by $n$ equations. Thus 
combining Proposition 9.12 (page 154) of \cite{F} with Example 9.1.5. 
we get (for a surface in ${\P}^{4}$):

\[
deg(c_{2}(T_{X}))+(\sigma_{1}(g_{i})-5)deg(c_{1}(T_{X}))+ \]
\beq
\label{inter}
+(\sigma_{2}(g_{i})
-5\sigma_{1}(g_{i})+15)deg(X)+W(p_{1},\ldots,p_{k})=\sigma_{4}(g_{i}),
\eeq
where $W(p_{1},\ldots,p_{k})$ is the weighted number of residual 
points and $\sigma_{j}(g_{i})$ is the j-th elementary symmetric 
polynomial in the degrees $g_{i}$ of the hypersurfaces $V_{i}$.
Assuming $X$ subcanonical, from $c_{1}(T_{X})=-K_{X}$, we get 
$c_{1}(T_{X})=-eH$, where $H$ is the class of a hyperplane section; 
moreover, from the exact sequence:
\[ 0\r T_{X}\r T_{{\P}^{4}}\otimes {\O} _{X} \r N_{X/{\P}^{4}} \r 0, \]
we get $c_{2}(T_{X})=10H^{2}+5HK+K^{2}-c_{2}(N_{X/{\P}^{4}})$, and by 
$ c_{1}(T_{X})=-eH$, we have 
$c_{2}(T_{X})=(10-5e+e^{2})H^{2}-c_{2}(N_{X/{\P}^{4}})$.
Since $deg(c_{2}(N_{X/{\P}^{4}}))=d^{2}$, $deg(H^{2})=d$ and 
$deg(H)=d$, (where $d$ is $deg(X)$), substituting in (\ref{inter}), 
taking degrees and  imposing $W(p_{1},\ldots,p_{k})=0$ reads:
\beq
\label{inter2}
[25+\sigma_{2}(g_{i})-(5+e)\sigma_{1}(g_{i})+e^{2}-d]d=\sigma_{4}(g_{i}).
\eeq
Thus, the relation (\ref{inter2}) gives a sufficient condition for a 
subcanonical surface in ${\P}^{4}$ to be scheme-theoretically defined 
by $4$ equations.

\subsection{An application: the linkage criterion}
As usual, if $X$ and $Y$ are l.c.i. of codimension $2$ in ${\P}^{n}$, 
we say that $X$ is (directly) linked to $Y$ if there exists 
a complete intersection $(F_{1},F_{2})$, such that $Y$ is the 
residual scheme of $X$ in the intersection $F_{1}\cap F_{2}$, and 
viceversa.
In \cite{BE}, working in characteristic zero and 
assuming $X$ smooth and subcanonical $dimX\geq 1$, Beorchia and Ellia proved that 
$X$ is a complete intersection iff it is self-linked, i.e. iff there 
exists complete intersection $(F_{1},F_{2})$ such that $F_{1}\cap 
F_{2}=2X$ ($F_{1}$ and $F_{2}$ define on $X$ a double 
structure which is a complete intersection). They also asked if the same criterion holds also for 
possibly singular l.c.i..  Recently, in \cite{FKL}, Franco, Kleiman 
and Lascu have given a positive answer to this question proving that the 
same criterion holds avoiding smoothness: $X$ can be reducible and 
non reduced. Their proof works only in 
characteristic zero (unless $dimX\geq 4$, where it holds over any 
algebraically closed field, due to a previous result of Faltings), so they ask if the same holds in positive 
characteristic, for lower dimensional $X$. Using our Theorem A we prove the following:

{\bf Proposition A}: {\em Let $X$ be a subcanonical (possibly singular) l.c.i. subscheme of 
codimension 2 in ${\P}^{n}$, $n\geq 4$, defined over an algebraically 
closed field of any characteristic. Then $X$ is a complete intersection iff it is 
self-linked.}

Proof: According to the ÒGherardelli linkage theoremÓ, which holds 
over
any algebraically closed field (see \cite{FKL} for its proof) we know that $X\subset 
F_{1}\cap F_{2}$ is subcanonical iff its residual scheme $Y$ (in the 
complete intersection $F_{1}\cap F_{2}$) is 
scheme-theoretically defined by the intersection of $F_{1}$ and 
$F_{2}$ with a third hypersurface $F_{3}$. On the other hand, if $X$ 
is self-linked by $F_{1}$ and $F_{2}$, then, by definition $X$ is 
equal to its own residual scheme in the complete intersection of 
$F_{1}$ and $F_{2}$, and since $X$ is assumed subcanonical, by the 
Gherardelli theorem it is scheme theoretically defined by $F_{1}, 
F_{2}$ and $F_{3}$; hence, by Theorem A, it is a complete 
intersection as soon as $dimX\geq 2$. Viceversa, if $X$ is a complete 
intersection, it is immediate to see that it is self-linked (just 
consider the intersection of $F_{1}$ and $2F_{2}$ if $X=F_{1}\cap F_{2}$). \vuoto

There is an immediate generalization of the previous proposition, 
which is the following:

{\bf Proposition B}: {\em Let $X$ as in Proposition A. Then $X$ is a 
complete intersection iff it can be (directly) linked to $Y$, where 
$Y$ is any subcanonical (possibly singular) l.c.i. subscheme.}

Proof: It is sufficient to use again the Gherardelli linkage and 
Theorem A. \vuoto

{\bf Remark 3}: The $X$'s as in Proposition B are self-linked iff 
they are scheme-theoretically defined by three hypersurfaces. Indeed, 
if $X$ is self-linked, then by Gherardelli it is ÒschematicallyÓ 
defined by $3$ equations; viceversa, if $X$ is defined by $3$ 
equations it is a complete intersection by Theorem A and then it is 
self-linked by Proposition B.

{\bf Remark 4}: The most difficult case, in order to characterize 
complete intersection via self-linking is that of curves in ${\P}^{3}$.
Beorchia and Ellia proved their criterion (in characteristic zero) also for
 curves, assuming that they are smooth, while Franco, Kleiman and 
 Lascu extended this result to l.c.i. curves (always working in 
 characteristic zero). In positive characteristic (characteristic 2),
  however, there is 
 certainly a counterexample for this criterion to hold in the case of 
 curves, due to Migliore (see 
 the discussion at the end of \cite{FKL}).
 So our extension of this criterion over a field of any characteristic 
 is the best possible for low dimensional subvarieties: that is 
 surfaces are the lowest dimensional subvarieties where this 
 criterion holds without exceptions.
 Unfortunately, up to now, there is no positive result for the case 
 of space curves in characteristic grater than zero.

\bigskip

{\bf Acknowledgements}: It is a pleasure to thank Philippe Ellia and 
Alexandru Lascu for useful discussions, and Steve Kleiman for 
valuable remarks and corrections; I would like 
also to thank the members of the Department of Mathematics of Ferrara University 
for the kind hospitality and the warm and stimulating atmosphere.

\end{document}